\documentclass[12pt]{article}
\usepackage{enumerate,amsthm,amsfonts,amssymb,pstricks,hyperref, topcapt}%
\usepackage[leqno]{amsmath}
\usepackage{amscd}
\usepackage{amsfonts,rotating}
\usepackage{graphicx}%
\usepackage{fancyhdr}

\divide\oddsidemargin
by 2  

\makeatletter

\newcommand{\address}[2][]{%
  \ifx\@add@ress\@undefined\gdef\@add@ress{\par\par\bigskip}\AtEndDocument{\@add@ress}\fi
  \g@addto@macro\@add@ress{\bigskip\noindent{\small\scshape%
      \ifx#1\empty\else{\bfseries Address of #1:}\ \fi#2}\par\par}}

\renewenvironment{abstract}{\small\quotation\noindent
  {\bfseries \abstractname}}{\endquotation \par}
\newcommand{\footnotetextplain}[1]{\begingroup\def\@thefnmark{}%
  \long\def\@makefntext##1{\parindent 0pt\noindent ##1}\@footnotetext{#1}
  \endgroup}

\newcommand{\AMSsubjclass}[2]{\footnotetextplain{2000
   \emph{Mathematics Subject Classification:} Primary #1, Secondary #2.}}

\newcommand{\keywords}[1]{\footnotetextplain{\emph{Key words and phrases:} #1.}}

\makeatletter

\xdef\qedbuit{\qed}

\newcommand{\TeoremaAmbFinalMarcat}[1]{%
  \expandafter\gdef\csname end#1\endcsname{\qedbuit\@endtheorem}}

\theoremstyle{definition}
 \TeoremaAmbFinalMarcat{defi}
 \TeoremaAmbFinalMarcat{ex}
 \TeoremaAmbFinalMarcat{rem}
 \TeoremaAmbFinalMarcat{hrem}
 \TeoremaAmbFinalMarcat{conv}
 \TeoremaAmbFinalMarcat{prob}
 \TeoremaAmbFinalMarcat{conj}
 \TeoremaAmbFinalMarcat{conj}

 \newenvironment{proclama not emphasized}[1]{\par\vspace{\topsep}\noindent{\bf #1}}{\par\vspace{\topsep}}

\newenvironment{prooftext}[1] 
               {\noindent \textit{ #1}~ }     
               {\hfill\rule{2.5mm}{2.5mm} \vspace{\parskip} } 

\newcommand{\start}[2]{\begin{#1}\label{#2}}
\newcommand{\secc}[1]{Section~\ref{#1}}
\newcommand{\theoc}[1]{Theorem~\ref{#1}}
\newcommand{\propc}[1]{Proposition~\ref{#1}}
\newcommand{\coryc}[1]{Corollary~\ref{#1}}
\newcommand{\lemc}[1]{Lemma~\ref{#1}}
\newcommand{\conjc}[1]{Conjecture~\ref{#1}}

\newcommand{\figc}[1]{Figure~\ref{#1}}

\newcommand{ \refc}[1]{~\ref{#1}}

\def\@enum@{\list{\csname label\@enumctr\endcsname}%
           {\usecounter{\@enumctr}\def\makelabel##1{\hss\llap{##1}}
           \itemsep=2pt\parsep=0pt\topsep=3pt plus 1pt minus 1 pt}}

\newenvironment{romlist}{\enumerate[(i)]}{\endenumerate}
\newenvironment{numlist}{\enumerate[(1)]}{\endenumerate}

\newcommand{\si}{\mathop\mathrm{SI}}

\newcommand{\A}{A}
\newcommand{\B}{B}
\newcommand{\ww}{\mathsf w}
\newcommand{\uu}{\mathsf u}
\newcommand{\vv}{\mathsf v}

\title{Self-intersection numbers of curves\\ in the
doubly-punctured plane}
\author{Moira Chas and  Anthony Phillips}
\date{\today}

\begin{document}
\maketitle

\begin{abstract} 
We address the problem of computing  bounds for the self-intersection
number (the minimum number of self-intersection points) of members 
of a free homotopy 
class of curves in the doubly-punctured plane as a function
of their combinatorial length $L$; this is
the number of 
letters required for a minimal description of the class in terms of the
standard 
generators of the fundamental group and their inverses. We prove that
the self-intersection
 number is bounded above by $L^2/4 + L/2 -1$, and that when $L$ is even,
this bound is sharp; in that case there are exactly four distinct 
classes attaining that bound.  When $L$ is odd, we establish
a smaller, conjectured upper bound ($(L^2-1)/4)$) in certain cases; and
there we show it is sharp. 
Furthermore, for the doubly-punctured plane, these self-intersection
numbers are
bounded  {\em below}, by $L/2-1$ if $L$ is even, $(L-1)/2$ if $L$ is odd;
these bounds are sharp.

\end{abstract}

\AMSsubjclass{57M05}{57N50,30F99}

\section{Introduction}

By the doubly-punctured plane we refer to the compact surface
with boundary (familiarly known as the ``pair of pants'')
obtained by removing, from a closed two-dimensional disc,
two disjoint open discs. 
This work extends to that surface
the research
reported in \cite{torus} for the punctured torus.  
Like the punctured torus, the doubly-punctured plane
has the homotopy type of a figure-eight.
Its fundamental group is free on two generators: once 
these are chosen, say  $a, b$, 
a free homotopy class of curves on the surface can
be uniquely 
represented as a reduced cyclic word in the symbols
$a, b, \A, \B$ (where $\A$ stands for $a^{-1}$ and $\B$ for $b^{-1}$). 
A {\em cyclic word } $w$ is an equivalence class of words 
related by a cyclic permutation of their letters; we will write
$w=\langle r_1 r_2 \dots r_n \rangle$ where the $r_i$ are
the letters of the word, and $\langle r_1 r_2 \dots r_n \rangle =
\langle r_2 \dots r_n r_1 \rangle$, etc. {\em Reduced}
means that the cyclic word contains no 
juxtapositions of $a$ with $A$, or $b$ with $B$.
The {\em length} (with respect to the generating set $(a,b)$) of 
a free homotopy class of curves is the number of letters occurring in
the corresponding reduced cyclic word.

This work studies the relation between length and
the {\em self-intersection
number} of a free homotopy class of curves: the 
smallest number of 
self-intersections among all general-position curves in the class. 
(General position in this context means as usual that there are
no tangencies or multiple intersections).
The self-intersection number is a property of the free homotopy class
and hence 
of the corresponding reduced cyclic word $w$; 
we denote it by $\si (w)$. 
 Note that a word and its inverse
have the same self-intersection number.

\begin{figure}[htp]
\centering \includegraphics[width=3in]{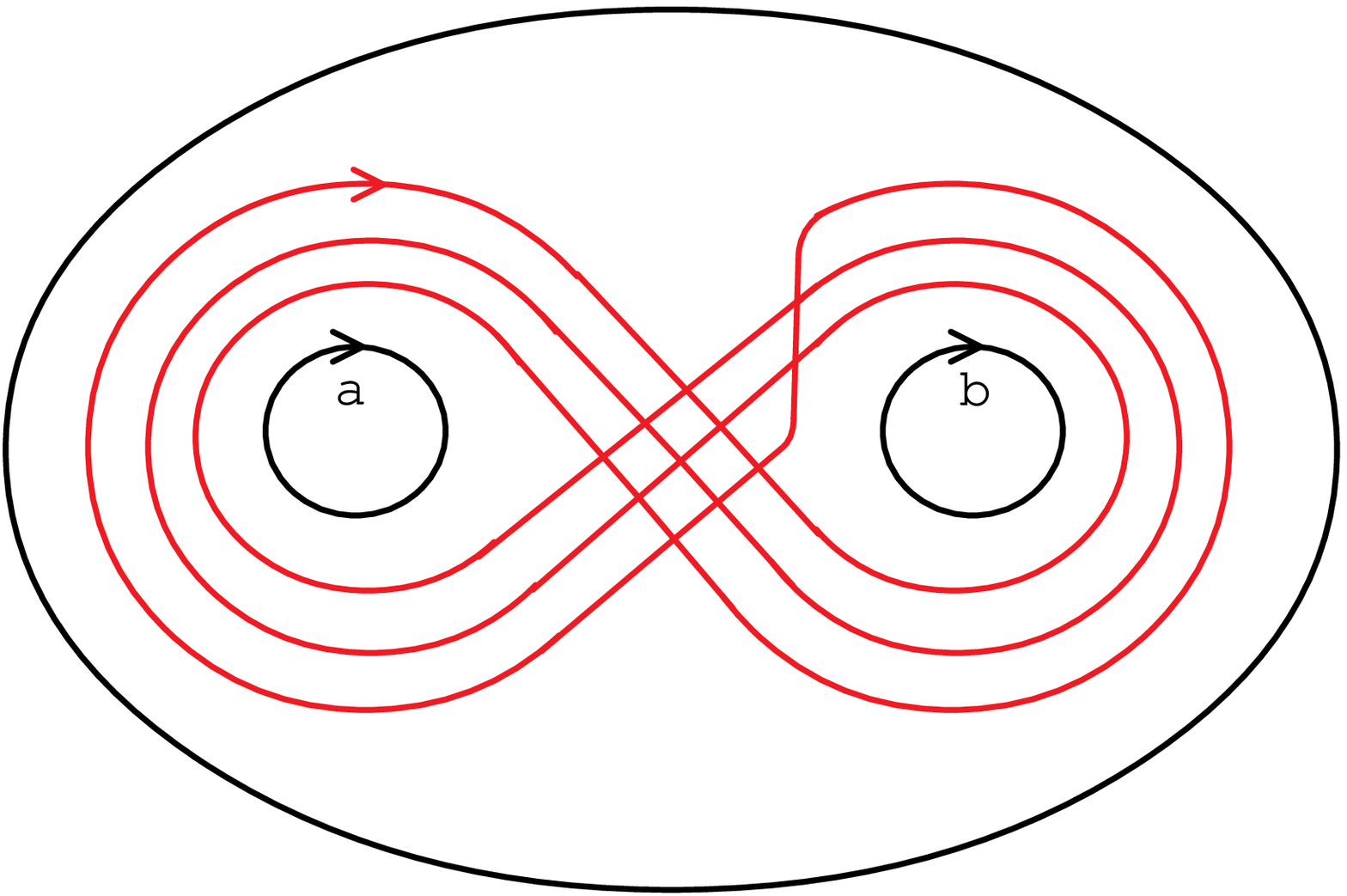} \includegraphics[width=3in]{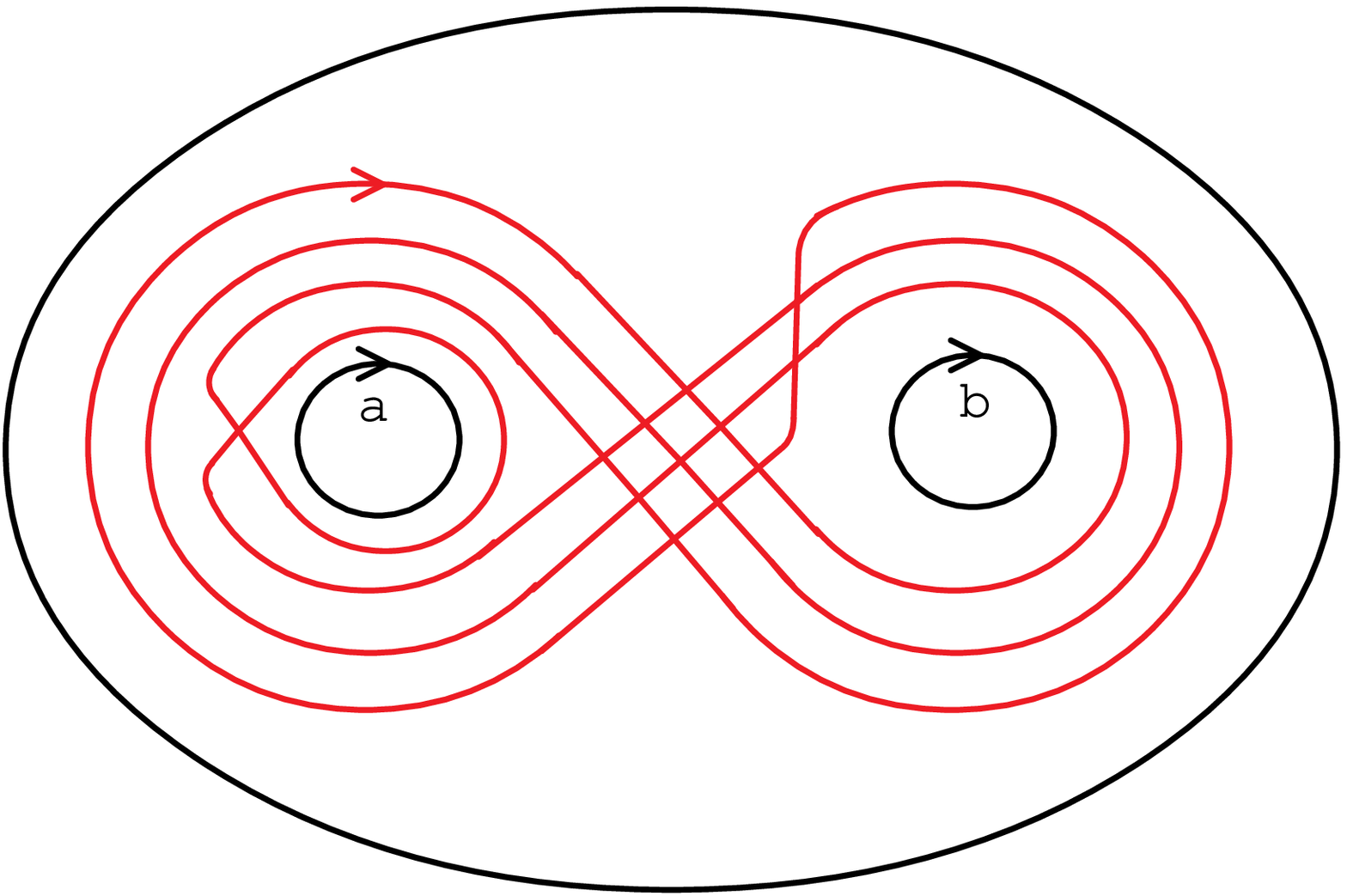}
\caption{Left: curves of the form $\langle aBaBaB\rangle$ 
have maximum self-intersection number $L^2/4 + L/2 -1$  
for their length  (\theoc{upper bound even}). Right:
curves of the form $\langle aaBaBaB\rangle$ 
have self-intersection number $(L^2 -1)/4$.  We
conjecture (\conjc{upper bound odd}) this is maximal, and
prove this conjecture in certain cases (\theoc{theo-odd}). }
\label{pantsgrids}
\end{figure}

\start{theo}{upper bound even}
\begin{numlist}
\item The self-intersection
number for a reduced cyclic word of 
length $L$ on the doubly-punctured plane is bounded above by
$L^2/4 + L/2 -1$. 
\item If $L$ is even, this bound is sharp: for  
$L \ge 4$ and even,  
the cyclic words realizing the maximal self-intersection number are 
(see \figc{pantsgrids}) $(aB)^{L/2}$ and $(Ab)^{L/2}$.
For $L=2$, they are $aa, AA, bb, BB, aB$ and $Ab$.
\item If $L$ is odd, the maximal self-intersection number of  words
of length $L$ is 
at least
 $(L^2-1)/4$.
\end{numlist}
\end{theo}

\start{conj}{upper bound odd} The maximal self-intersection
number for a reduced cyclic word of 
odd length $L = 2k+1$ on the doubly-punctured plane is
$(L^2-1)/4$; the words realizing the maximum have
one of the four forms $\langle (aB)^kB\rangle, \langle a(aB)^k\rangle, 
\langle (Ab)^kb\rangle, \langle A(Ab)^k\rangle.$
\end{conj}

\start{defi}{blocks} Any reduced cyclic word is either a pure power
or may be written in the form
$\langle\alpha_1^{a_1}\beta_1^{b_1}\dots\alpha_n^{a_n}\beta_n^{b_n}\rangle$, 
where
$\alpha_i \in \{a, A\}$, $\beta_i \in \{b, B\}$ , all $a_i$ and $b_i$
are positive, and 
$\sum_1^n(a_i + b_i) = L$, the length of the word.   We will
refer to each $\alpha_i^{a_i}\beta_i^{b_i}$ as an {\em $\alpha\beta$-block},
and to $n$ as the word's {\em number of $\alpha\beta$-blocks}.
\end{defi}

\start{theo}{theo-odd} On the doubly-punctured plane, consider 
 a reduced cyclic word $w$ of odd length $L $ with $n$
$\alpha\beta$-blocks. If $L>3n$, or $n$ is prime, or $n$ is a power of $2$, then
the self-intersection number of $w$ 
satisfies $\si(w)\leq \frac{L^2 -1}{4}$. This bound is sharp.
\end{theo}

The doubly punctured plane has the property that self-intersection
numbers of words are bounded {\it below}.

\start{theo}{lower-bound} On the doubly punctured plane, curves in
the free homotopy class  represented by
a reduced cyclic word of length $L$ have at least $L/2-1$ self-intersections if $L$ is even and $(L-1)/2$ self-intersections if $L$ is odd. These bounds are achieved by $(ab)^{\frac{L}{2}}$ and $(AB)^{\frac{L}{2}}$ if $L$ is even and by the four words $a(ab)^{\frac{L-1}{2}}$, etc. when $L$ is odd.

\end{theo}

\start{cory}{finite} A curve with minimal self-intersection number $k$ has
combinatorial length at most $2k+2$. There are therefore only 
finitely many free homotopy classes with 
minimal self-intersection number $k$. 
\end{cory}

\start{rem}{unique}   A surface of negative Euler characteristic which is not
the doubly punctured plane has infinitely many homotopy classes of simple
closed curves \cite{mm}. Since the $(k+1)$st power of a simple closed curve 
kas self-intersection number $k$, it follows that for any $k$ there
are infinitely many distinct homotopy classes of curves with self-intersection
number $k$. (A   more elaborate argument using the 
mapping class group constructs, for any $k$, infinitely many 
distinct {\em primitive}
classes (not a proper power of another class) with self-intersection 
number $k$). So the doubly punctured plane is the unique surface of
negative Euler characteristic satisfying \coryc{finite}.
\end{rem}

The authors have benefited from discussions with Dennis Sullivan, and
are very grateful to Igor Rivin who contributed an essential
element to the 
proof of \theoc{lower-bound}. Additionally, they have profited from
use of Chris Arettines' JAVA program, which draws minimally self-intersecting
representatives of free homotopy classes of curves in surfaces.
The program is currently available at

http://www.math.sunysb.edu/$\sim$moira/applets/chrisApplet.html

\subsection{Questions and related results}\label{rr}

The doubly punctured plane admits a 
hyperbolic metric making its boundary geodesic. An elementary argument 
shows that for curves on that surface, hyperbolic and combinatorial
lengths are quasi-isometric. Some of our combinatorial results can
be related  in this way to statements about intersection numbers and
 hyperbolic length.

 A free
homotopy class of combinatorial length $L$ in a  surface with boundary can be represented by $L$ chords in a fundamental polygon. Hence,
 the maximal self-intersection number of a cyclic reduced word of length $L$ is bounded above by $\frac{L(L-1)}{2}$.

We prove in \cite{torus} that for the punctured
torus the maximal self-intersection number  $\si_{\max}(L)$ of a free homotopy
class of combinatorial length $L$  is equal to  $(L^2-1)/4$ if $L$ is 
even and to $(L-1)(L- 3)/4$ if $L$ is odd. 
This implies that the limit of $\si_{\max}(L)/L^{2}$ is $\frac{1}{4}$ as $L$ 
approaches infinity. (Compare \cite{lalley}). The same limit holds for the doubly punctured plane (\theoc{upper bound even}). 
On the other hand, according to our (limited) experiments, there are no analogous polynomials for more  
general surfaces; but it seems reasonable to ask:

\start{ques}{xxx} 
Consider closed curves on a surface $S$ with boundary. 
Let $\si_{\max}(L)$ be the maximum self-intersection number for 
all curves of combinatorial
length $L$. Does  $\si_{\max}(L)/L^{2}$ converge? And if so, to what limit?
Does this limit approach $\frac{1}{2}$ as the genus of the surface
approaches infinity? \end{ques}

\start{ques}{xxy}
Consider closed curves on a hyperbolic surface $S$ (possibly closed). 
Let $\si_{\max}(\ell)$ be the maximum self-intersection number for any
 curve of 
{\em hyperbolic}
length at most $\ell$. Does  $\si_{\max}(\ell)/\ell^{2}$ converge? And if so, to what limit?

\end{ques}

Basmajian \cite{basmajian} proved for a closed, hyperbolic surface $S$
that  there exists an increasing sequence $M_k$ (for $k = 1, 2, 3, ...$) 
 going to
infinity so that if $w$ is a closed geodesic with self-intersection number 
$k$, then its geometric length is larger than $M_k$ . Thus the length of a closed geodesic gets arbitrarily large 
as its self-intersection gets large. For the
doubly punctured plane, in terms of the combinatorial length, we calculate $M_{k}=\sqrt{5+4k}-1$.

\section{A linear model}\label{linear}

In this section we will need to distinguish between
a cyclically reduced linear word $\ww$ in the generators and their inverses,
and the the associated reduced cyclic word $w$. 
We 
introduce an algorithm for constructing 
from $\ww$
a representative curve  for  $w$.
An upper bound for the
self-intersection numbers of these representatives may be
easily estimated; taking the minimum of this bound
over cyclic permutations of $\alpha\beta$-blocks
 will yield 
  a useful upper bound for $\si(w)$.

\subsection{Skeleton words}\label{skel}

Given a 
cyclically reduced  word  
$w = \langle\alpha_1^{a_1}\beta_1^{b_1}\dots\alpha_n^{a_n}\beta_n^{b_n}\rangle$, 
where $\alpha_i = a \mbox{~or~} A$, $\beta_i = b \mbox{~or~} B$ 
and all $a_i, b_i > 0$,
the corresponding {\em skeleton word} is 
$w_S = \langle\alpha_1\beta_1\dots\alpha_n\beta_n \rangle$, a word of length $2n$.

We now describe a systematic way for drawing a representative curve
for  $w_S$ starting from one of its linear forms $\ww_S$, 
and for {\em thickening} this curve to a representative for $w$.

\begin{figure}[htb]
\centering \includegraphics[width =3in]{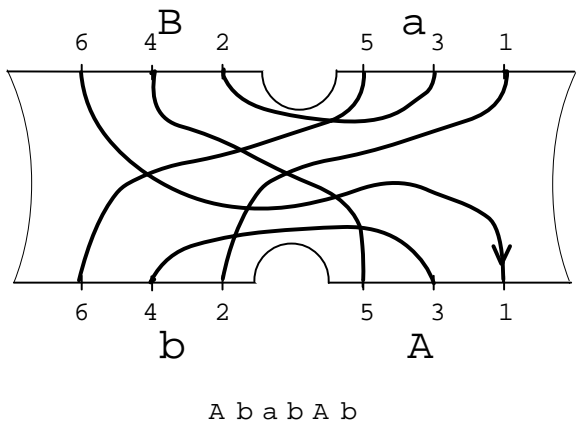}
\caption{The skeleton curve $AbabAb$. }
\label{skeleton1}       
\end{figure}

\begin{figure}[htp]
\centering
\includegraphics[width=3in]{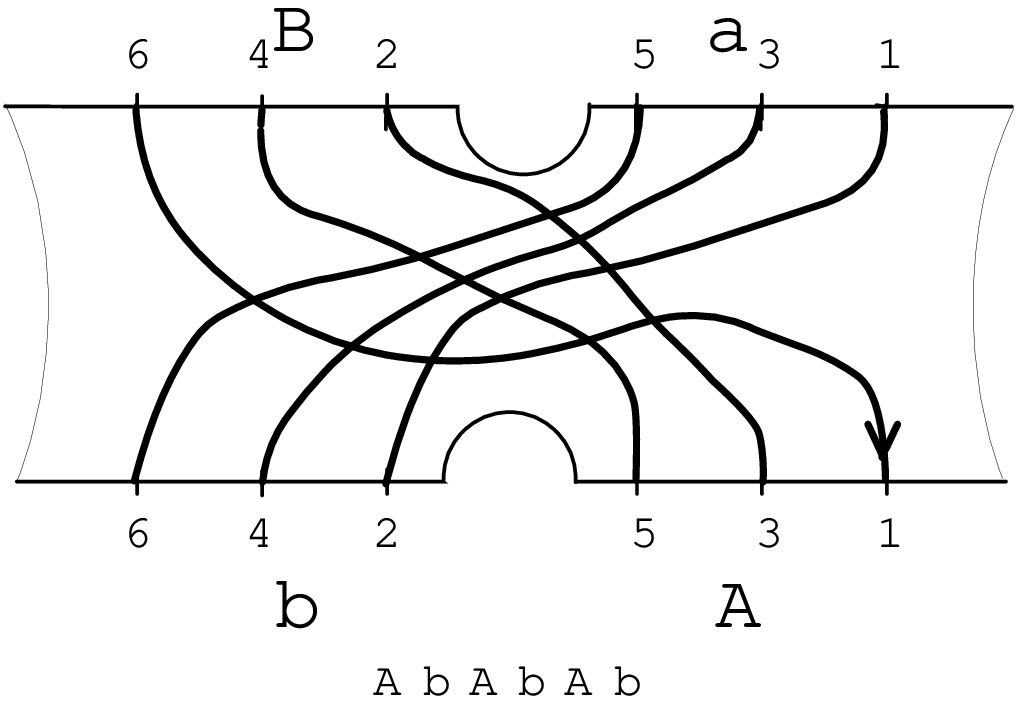}
\includegraphics[width=3in]{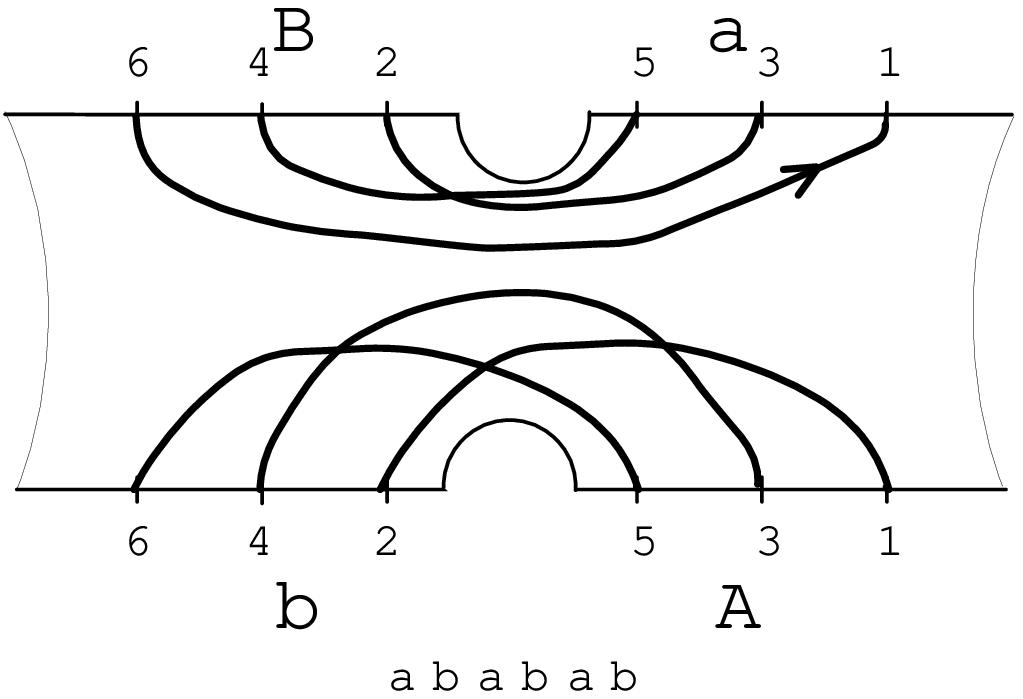}
\caption{The skeleton curves $ababab$ and $AbAbAb$. }
\label{skeleton-ababab}
\end{figure}

\begin{proclama not emphasized}{The 
skeleton-construction algorithm:}  (See Figures \ref{skeleton1} and \ref{skeleton-ababab}) Start by marking off $n$ points along each of the edges of the fundamental
domain; corresponding points on the $a,A$ sides are numbered
$1, 3, 5, \dots, 2n-1$ starting from their common corner; and
similarly corresponding points on the $b,B$ sides are
numbered $2n, \dots, 6, 4, 2$, the numbers decreasing away from
the common corner.

If the first letter  in 
$\ww_S$ is $a$, draw a curve segment entering
the $a$-side at 1, and one exiting the $A$-side at 1 (vice-versa
if the first letter is $A$). That segment is
then extended to enter the $b$-side at 2 and exit the $B$-side at 2
if the next letter in $\ww_S$ is $b$; vice-versa if it is $B$. And so
forth until the curve segment exiting the $b$ (or $B$)-side at $2n$
joins up with the initial curve segment drawn.
\end{proclama not emphasized}

We will refer to a segment of type $ab, ba, AB, BA$ as a {\em corner
segment}, and one of type $aB, Ab, bA, Ba$ as a {\em transversal}.
Note that (as above) a skeleton word has even length $2n$ and
therefore has $2n$ segments (counting the {\em bridging segment}
made up of the last letter and the first). 
The number of transversals must also be even, since if they are
counted consecutively they go from lower-case to upper-case or
vice-versa, and the sequence (upper, lower, ... ) must end up
where it starts. It follows that the number of corners is also
even.

\start{prop}{Abn} The self-intersection number of the representative of 
$(Ab)^{n}$ or $(aB)^n$ given by the curve-construction algorithm  
equals $n^2 + n - 1$.
\end{prop}

\begin{proof} Consider $(Ab)^{n}$; see \figc{skeleton-ababab}, left. 
This curve has only transversals.
There are $n$ parallel segments of type $Ab$; 
they join $1, 3, \dots, (2n-1)$ on the $a$-side
to $2, 4, \dots, 2n$ on the $b$-side. There are
$n-1$ parallel segments of type $bA$, which
join $2, 4, \dots, 2n-2$ on the $B$-side to
$3, 5, \dots 2n-1$ on the $A$-side. Each of these
intersects all $n$ of the $Ab$ segments. Finally
the bridging $bA$ segment joins
$2n$ on the $B$-side to $1$ on the $A$-side.
This segment begins to the left of all the
other segments and ends up on their right:
it intersects all $2n-1$ of them. The total
number of intersections is $n(n-1) + 2n -1 =
n^2 + n - 1$. A symmetrical argument handles $(aB)^n$.
\end{proof}

\start{prop}{abn} The self-intersection number of the representative of $(ab)^{n}$ given by the curve-construction  algorithm  equals $(n-1)^{2}$.\end{prop}


\begin{proof} (See \figc{skeleton-ababab}, right)
This curve has only corners. There are
$n$ segments of type $ab$, joining $1, 3, \dots, 2n-1$
on the $A$-side to $2, 4, \dots, 2n$ on the $b$-side.
Since their endpoints interleave, each of these curves
intersects all the others. There are $n-1$ segments of
type $ba$, joining $2, 4, \dots, 2n-2$ on the $B$-side
to $3, 5, \dots 2n-1$ on the $a$-side. Again, each of
these curves intersects all the others. Finally the
bridging $ba$ segment joining $2n$ to $1$ spans both endpoints
of all the others and so intersects none of them.
The total number of intersections is $\frac{1}{2}n(n-1)
+ \frac{1}{2}(n-1)(n-2) = (n-1)^2$. 
\end{proof}

\start{prop}{prop-2c} Let 
$w$
be a skeleton word of length $2n$. 
The number of corner segments
in $w$ is even, as remarked above;
write it as $2c$.
Then the self-intersection
number of 
$w$
 is bounded above by $n^2 + n - 1 - 2c$.
\end{prop}

\begin{proof}
Using Propositions \ref{Abn} and \refc{abn} we can assume that 
$w$ has both corner-segments and  transversals.
We may then choose a linear representative $\ww$ 
with the property 
that the
bridging segment between the end of the word and the
beginning is a transversal.
Of the $2c$ corners, $c$ will be on top (those of type $AB$ or $ba$) 
and $c$ on the bottom (types $ab$ and $BA$). An $ab$ or $AB$ corner segment
joins a point numbered $2j-1$ to a point numbered $2j$ on the same
side, top or bottom, as $2j-1$. 
It encloses segment endpoints $2j+1, 2j+3, \dots, 2n-1, 2, 4, 
\dots 2j-2$, a total of $n-1$ endpoints; 
similarly, a $ba$ or $BA$ segment
encloses $n-2$ endpoints. So there are at most $2c(n-1) -c(c-1)$ intersections
involving corners, correcting for same-side corners
having been counted twice. The $2n-2c$ transversals intersect each
other just as in the pure-transversal case, producing $(n-c)^2 + (n-c) -1$
intersections. The total number of intersections is therefore
bounded by $n^2 + n -1 -2c$. \figc{skeleton1} shows
the curve $AbabAb$ (here $n=3, c=1$) with 8 self-intersections.
\end{proof}

\subsection{Thickening a skeleton; proof of \theoc{upper bound even} (1), (2)}
Once the skeleton curve corresponding to $\ww_S$ is constructed, 
it may be {\em thickened} to
produce a representative curve for $w$. The algorithm
runs as follows. 

\begin{figure}[htp]
\centering
\includegraphics[width=5.5in]{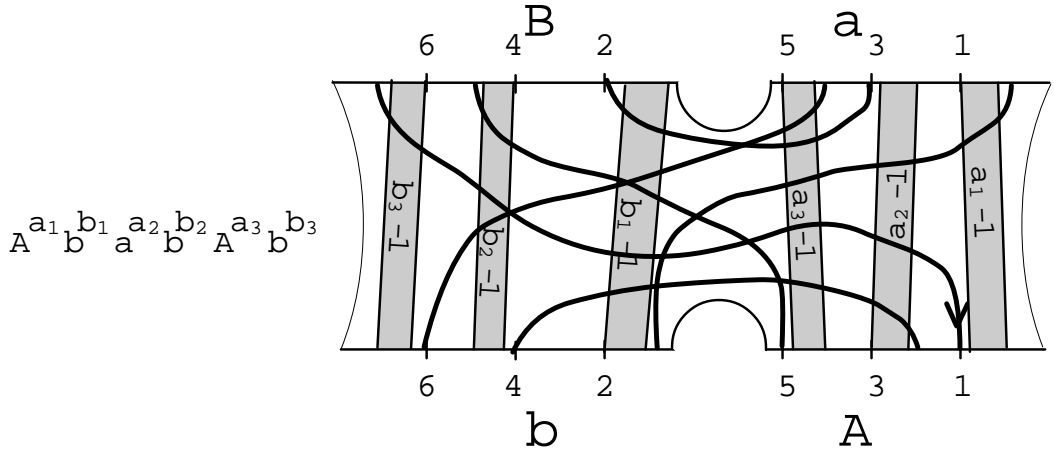}
                         

\caption{The skeleton curve $AbabAb$ thickened to
represent the linear word $A^{a_1}b^{b_1}a^{a_2}b^{b_2}A^{a_3}b^{b_3}$. 
The grey bands represent the curve segments corresponding to the
extra letters: $a_1-1$ copies of $A$, etc. Notice that the segments
from the skeleton curve intersect the $a_1$ and $b_3$ bands once,
the $a_2$ and $b_2$ bands three times, and the $a_3$ and $b_1$
bands five times.}\label{thickening}     
\end{figure}

\begin{proclama not emphasized}{The skeleton-thickening algorithm.} (See \figc{thickening}) Suppose for explicitness that $w$ starts with $A^{a_1}$. The extra $a_1-1$
copies of $A$, inserted after the first one, correspond to segments
entering the $a$-side (the first one at 1)  and exiting the $A$-side
(the last one at a point opposite the displaced entrance point
of the first skeleton segment); the new segments are parallel.
Similarly the extra $b_1-1$ segments appear as parallel segments
originating and ending near the 2 marks on the $b$ and $B$-sides;
so there are no intersections between these segments and those in
the first band. Proceeding in this manner we introduce $n$ non-
intersecting bands
of $a_1-1, b_1-1, a_2-1, ..., b_n-1$ parallel segments. New intersections
occur between these bands and segments of the skeleton curve. The
two outmost bands (corresponding to $a_1$ and $b_n$) are each 
intersected by one of the skeleton segments; the next inner bands
($a_2$ and $b_{n-1}$) each intersect three of the skeleton
segments; \dots;  the two innermost bands ($a_n$ and $b_1$)
each intersect $(2n-1)$
of the skeleton segments. 
\end{proclama not emphasized}

Adding these intersections to the
bound on the self-intersections of the skeleton curve itself yields
$$\si(w) \leq (a_1 + b_n - 2) + 3(a_2 + b_{n-1} - 2) + \cdots +
(2n-1)(a_n + b_1 - 2) + n^2 + n - 1.$$  
Since $1 + 3 + \cdots + (2n-1) = n^2$ we may repackage this
expression as
$$\si(w) \leq f(a_1, \dots, a_n, b_1, \dots, b_n) - n^2 + n - 1,$$
where we define $f$ by
$$f(a_1, \dots, a_n, b_1, \dots, b_n) = (a_1 + b_n) + 3(a_2 + b_{n-1})
+ \cdots + (2n-1)(a_n + b_1).$$
\vspace{.1in}

Applying the skeleton-thickening algorithm to
the cyclic permutation $\alpha_1^{a_1}\beta_1^{b_1}\dots\alpha_n^{a_n}\beta_n^{b_n} \rightarrow 
\alpha_2^{a_2}\beta_2^{b_2}
\dots\alpha_n^{a_n}\beta_n^{b_n}\alpha_1^{a_1}\beta_1^{b_1}$ yields another
curve representing the same word. There are $n$ such permutations, leading to
\begin{equation}\label{equation(*)}
\si(w) \leq [\min_{i=0,\dots,n-1} 
f\circ r^i(a_1, \dots, a_n, b_1, \dots, b_n)] - n^2 + n - 1,
\end{equation}

where $r$ is the coordinate permutation
$(a_1, \dots, a_n, b_1, \dots, b_n)  \rightarrow
(a_2, \dots, a_n, a_1, b_2, \dots, b_n, b_1).$

\start{prop}{nL} Set $L = a_1 + \cdots + a_n + b_1 + \cdots + b_n$.
Then ${\displaystyle \min_{i=0,\dots,n-1}
f\circ r^i(a_1, \dots, a_n, b_1, \dots, b_n) \leq nL}.$
\end{prop}

\begin{proof}
We write
$$f(a_1, \dots, b_n) = (a_1 + b_n) + 3(a_2 + b_{n-1})+ \cdots + (2n-1)(a_n + b_1)$$
$$f\circ r(a_1, \dots, b_n) = (a_2 + b_1) + 3(a_3 + b_n)+ \cdots + (2n-1)(a_1 + b_2)$$
$$.$$
$$.$$
$$.$$
$$f\circ r^{n-1}(a_1, \dots, b_n) = (a_n + b_{n-1}) + 3(a_1 + b_{n-2})+ \cdots + (2n-1)(a_{n-1} + b_n).$$

The average of these $n$ functions is
$\frac{1}{n}(L + 3 L + \cdots (2n-1)L) = nL.$
Since the minimum of $n$ functions must be less than their average,
the proposition follows.
\end{proof}

\begin{prooftext}{Proof of \theoc{upper bound even}, (1) and (2)}
We work with 
$w = 
\langle\alpha_1^{a_1}\beta_1^{b_1}\dots\alpha_n^{a_n}\beta_n^{b_n}\rangle$. 
We have established that
$$\si(w) \leq \min_{i=0,\dots,n-1} 
f\circ r^i(a_1, \dots, a_n, b_1, \dots, b_n) - n^2 + n - 1.$$
Using Proposition \ref{nL},
$$\si(w) \leq nL - n^2 + n - 1 = -n^2 + n(L+1) -1.$$
For a given $L$, this function has its real maximum at $n = (L+1)/2$.
Since each 
$\alpha\beta$-block contains at least 2 letters, $n$ must be less than or 
equal to $L/2$. So a
bound on $\si(w)$ is the value at $n = L/2$ ($L$ even) or 
$n = (L - 1)/2$ ($L$ odd):
$$\si(w) \leq \left \{ 
\begin{array}{ll}
L^2/4 + L/2 -1 & (L \mbox{ even})\\
L^2/4 + L/2 -7/4 & (L \mbox{ odd}).
\end{array} \right .  $$ 

For $L$ even, note (\propc{Abn})
that the skeleton words  $w=(aB)^n$ and $w=(Ab)^n$ 
satisfy  $\si(w) = n^2 + n - 1 = L^2/4 + L/2 -1$; so the 
bound for this 
case is sharp; furthermore since words with $n=L/2$ must be
skeleton words, it follows from \propc{prop-2c} 
these are the only words attaining the bound.

\end{prooftext}

\start{rem}{odd}
 For $L$ odd, our numerical experiments (which go up to $L=20$) and the special cases we
prove below have $\si(w) \leq (L^2-1)/4 $, so the function
constructed here does not give a sharp bound.

\end{rem}

\section{Odd length words}\label{odd-length}

\subsection{A lower bound for the maximal  self-intersection number;
proof of \theoc{upper bound even} (3)}

\begin{prooftext}{Proof of \theoc{upper bound even}, (3) 
(The maximum self-intersection number
for words of odd length $L$ is at least $(L^2-1)/4$).}
We will show that the words of the form $a(aB)^{\frac{L-1}{2}}$ have self-intersection equal to $(L^{2}-1)/4$. Consider a 
representative of $w$ as in \figc{pantsgrid2}, 
where $n = \frac{L-1}{2}$. There is 
an $n \times n$ grid of intersection points in the center,
plus the $n$ additional intersections $p_2, \dots p_{2n}$, a total
of $n^2 + n = (L^{2}-1)/4$. We need to check that none of these 
intersections spans a bigon (this is the only way
\cite{hs} that an intersection can be deformed away).

\begin{figure}[htp]
\centering
\includegraphics[width=2.5in]{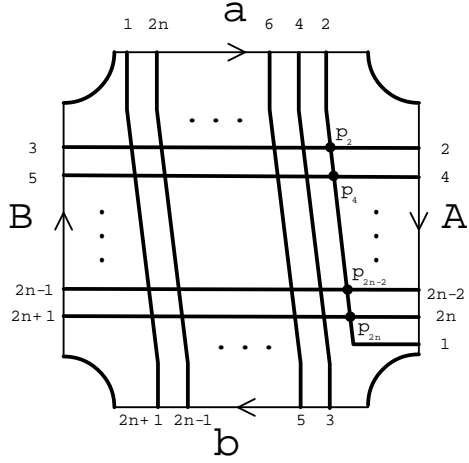}

\caption{The curve $a(aB)^n$ 
represented in the fundamental domain for the doubly punctured
disc.}\label{pantsgrid2}     
\end{figure}

With notation from \figc{pantsgrid2}, the only vertices that could be
part of a bigon are those from which two segments exit along the same
edge, i.e. $p_2, p_4, \dots, p_{2n}$ . If we follow the segments from
$p_2$ through edge $A$ they lead to 1 on edge $A$ and $2n+1$ on edge
$b$, so no bigon there; the segments from $p_4$ through edge $A$
lead to $3, 2n+1$ on edge $b$, to $2, 2n$ on edge $A$ and then to
1 on edge $A$ and $2n-1$ on edge $b$, so no bigon; etc. Finally
the segments from $p_{2n}$ through edge $A$ lead to $2n-1, 2n+1$
on edge $b$ and eventually to 1 on edge $A$ and 3 on edge $b$:
no bigon.   
\end{prooftext}

\subsection{Preliminaries for upper-bound calculation}\label{odd-prelim}
In the analysis of self-intersections of odd length
words the exact relation between $L$ (the length
of a word) and $n$  (its number of $\alpha\beta$-blocks)
becomes more important. 

\start{prop}{3n bound} If a word $w$ has length $L$ and
$n$ $\alpha\beta$-blocks, with  $L \geq 3n$, then 
$\si(w) \leq \frac{1}{4}(L^2 -1).$ 
Note that by \theoc{upper bound even} (3),
this estimate is sharp.
\end{prop}
\begin{proof}
As established in the previous section (equation \ref{equation(*)})
$\si(w) \leq nL - n^2 + n -1.$

The inequality $nL - n^2 + n -1 \leq \frac{1}{4}(L^2 -1)$ is
equivalent to $L^2 -4nL + 4n^2 - 4n +3 \geq 0$. As a function
of $L$ this expression has two roots: $ 2n \pm \sqrt{4n-3}$;
as soon as $L$ is past the positive root, the inequality is
satisfied.

If $n \geq 3$, then $ L\geq 3n$
implies $L \geq 2n + \sqrt{4n-3}$.

If $n = 2$ our inequality
$\si(w) \leq nL-n^2 + n -1$ translates to $\si(w) \leq 2L-3$
which is less than $\frac{1}{4}(L^2 -1)$ always. 

If $n=1$ our inequality becomes
 $\si(w) \leq L-1$, which is less than $\frac{1}{4}(L^2 -1)$
as soon as $L \geq 3$. The other only possibility is $L=2$, an
even length.
\end{proof}

\subsection{The cases:  $n$ prime or $n$ a power of 2; proof of
\theoc{theo-odd} }

Other results for odd-length words require a more detailed
analysis of the functions $f\circ r^i(a_1, \dots, a_n, b_1, \dots, b_n)$,
keeping the notation of the previous section.

The proof of the following results is straightforward.

\start{lem}{t_i}

For a fixed $(a_1, \dots, a_n, b_1, \dots, b_n)$, set 
$$s_a = a_1 + \cdots + a_n,$$
$$ s_b = b_1 + \cdots + b_n,$$
$$t_i = f\circ r^i(a_1, \dots, a_n, b_1, \dots, b_n).$$
 Then
\begin{romlist}
\item  $t_{i+1}-t_i = 2n(a_i-b_i) -  2(s_a-s_b)$
\item $t_0 - t_{n-1} = 2n(a_n-b_n) - 2(s_a-s_b)$.
\item $t_{i+j} - t_i = 2n(a_i+\cdots +a_{i+j-1} -
b_i-\cdots -b_{i+j-1})-
2j(s_a-s_b).$
\end{romlist}
In particular, if $t_i = t_{i+r}$, for some $r>0$, then
$$n(a_1-b_1+a_2-b_2+\cdots +a_{i+r-1}-b_{i+r-1}) = r (s_a - s_b).$$
\end{lem}

\start{lem}{diff ti} If $n$ is prime and $L < 3n$, then all
the numbers $t_0, \dots, t_{n-1}$ are different.  
\end{lem}

\begin{proof} By \lemc{t_i}, if $t_i = t_{i+r}$, for some $r>0$,
then $n$ must divide $r$ or $s_a - s_b$. We will show each is
impossible. The first cannot happen because $r < n$. As for the
second, observe that
$s_a \geq n$ and $s_b \geq n$, and that  their sum is $L < 3n$;
so $s_a - s_b = s_a + s_b - 2s_b < 3n - 2n = n$. So $n$ cannot divide
$s_a - s_b$ either.
\end{proof}

\start{lem}{diff ti-2} If $n$ is a power of $2$ and $L$ is odd,
then all
the numbers $t_0, \dots, t_{n-1}$ are different.
\end{lem}

\begin{proof} Arguing as in \lemc{diff ti}: 
in this case, since $r < n$ it cannot be a multiple
of $n$, so $s_a - s_b$ must be even. But $s_a - s_b$ is congruent
mod 2 to $s_a + s_b = L$, which is odd.
\end{proof}

\start{prop}{n} If a word $w$ of odd length $L$ has a number  
of $\alpha\beta$-blocks which is prime or a power of two then $\si(w)\le (L^{2}-2)/4$.
\end{prop}

\begin{proof}
Let $n$ be the number of $\alpha\beta$-blocks in $w$.
 By Lemmas \ref{diff ti} and \refc{diff ti-2} the numbers $t_0, \dots, t_{n-1}$ are all different; in
fact (\lemc{t_i}) their differences are all even, so any two of
them must be at least 2 units apart.  It follows that 
$$\sum_{i=0}^{n-1}t_i \geq \min t_i + (\min t_i + 2) + \cdots
+ (\min t_i + 2n-2) = n\min t_i + n(n-1)$$
so their average, which we calculated 
in the proof of \propc{nL} to be $nL$, is greater than or
equal to $\min t_i + n -1$, and so 
(using equation \ref{equation(*)})
$$\si(w) \leq \min t_i -n^2 + n - 1 \leq nL -n^2 = n(L-n) \leq L^2/4;$$    
since $L$ is odd and $\si(w)$ is an integer, this means
$$\si(w) \leq (L^2-1)/4.$$
\end{proof}

Propositions \ref{3n bound} and  \ref{n} prove \theoc{theo-odd}.

\section{Lower bounds; proof of \theoc{lower-bound}}\label{lower}

\start{defi}{positive} A word in the generators of a surface group and their inverses is \emph{positive} if no generator occurs along with its inverse. 
Note that a positive word is automatically cyclically reduced.
\end{defi}

\start{nota}{alpha1} If $w$ is a word in the alphabet $\{a, A, b, B\}$,
we denote by $\alpha(w)$ (resp. $\beta(w)$) the total number of occurrences 
of  $a$ and $A$ (resp. $b$ and $B$).
\end{nota}

\start{prop}{lower-case} For any reduced cyclic word $w$ in the
alphabet $\{a, A, b, B\}$ there is a positive
cyclic word $w'$ of the same length
with $\alpha(w') = \alpha(w), ~\beta(w') = \beta(w)$
and  $\si(w') \leq \si(w)$.
\end{prop}
\begin{proof} 
We show how to change $w$ into a word written with
only  $a$ and  $b$
while controlling the self-intersection number.
If all the letters in $w$ are capitals, take $w' = w^{-1}$.
Otherwise, look in $w$ for a maximal (cyclically)
connected string of (one or more) capital letters.
The letters at the ends of this string must be one of the pairs
$(A,A), (A,B), (B,A), (B,B)$. In the case $(B,B)$
(the other three cases admit a similar analysis), focus on that
string and write
 $$w=\langle xa^{a_1}B^{b_1}A^{a_2}B^{b_2}\dots A^{a_{i}}B^{b_{i}}a^{a_{i+1}}
\rangle$$ 
where $x$ stands for the rest of the word.

 Consider a representative of $w$ with minimal self-intersection. 
In this representative consider the arcs corresponding to the segments 
$aB$ (joining the last $a$ of the $a^{a_{1}}$-block to the first
 $B$ of $B^{b_{1}}$) and $Ba$ (joining the last $B$ in $B^{b_{i}}$ to the first $a$ in $a^{a_{i+1}}$). These two arcs intersect in a point $p$.  Perform surgery around $p$ in the following way: remove these two segments, and replace them with an $ab$ and a $ba$ respectively, using the same endpoints. This
 surgery links the arc $a^{a_{i+1}}xa^{a_{1}}$ to the arc
  $B^{b_1}A^{a_2}B^{b_2}\dots A^{a_{i}}B^{b_{i}}$ traversed in the 
opposite direction, i.e. gives
  a curve corresponding to the word 
$$ w' = \langle a^{a_{i+1}}xa^{a_{1}}
  (B^{b_1}A^{a_2}B^{b_2}\dots A^{a_{i}}B^{b_{i}})^{-1}\rangle. $$ 
This word has the same $\alpha$ and $\beta$
values as $w$, has lost at least one self-intersection, and has 
strictly fewer upper-case letters than $w$. The process may be 
repeated until all upper-case letters have been eliminated.
\end{proof}

\start{prop}{cut-the-word}
In any surface $S$ with boundary,  Let $w$ be a cyclically reduced word in 
the generators of $\pi_1S$ which does not admit a  simple representative curve.
Then a linear word $\ww$ 
representing $w$ (notation from \secc{linear})  
can be written 
as the concatenation $\ww=\uu\cdot \vv$ 
of two linear words, in such a way that the associated cyclic words
satisfy
$\si(u)+\si(v)+1 \leq \si(w)$. (Note that $u$ and $v$ are not necessarily cyclically reduced).
 \end{prop}

\begin{proof} 

\begin{figure}[htp]
\centering
\includegraphics[width=5.5in]{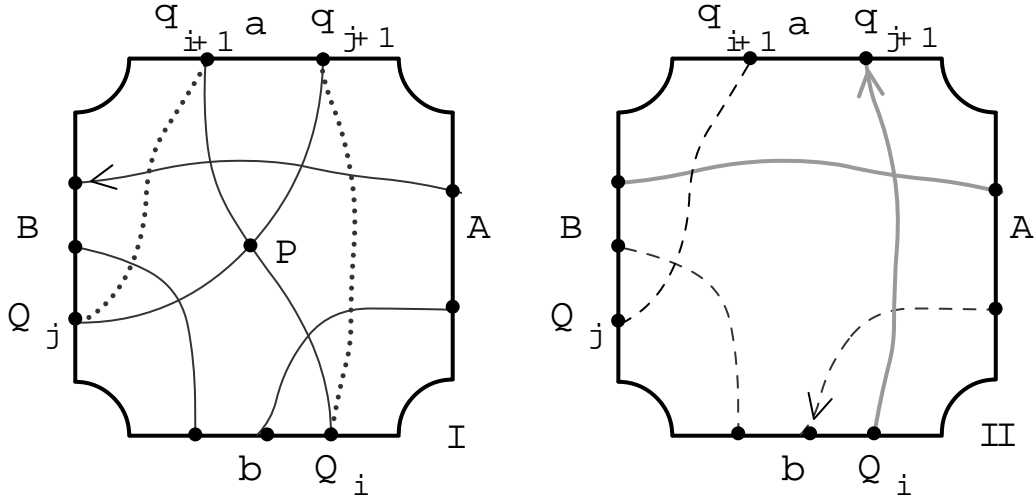}
\caption{Splitting $\ww$ as $\uu\cdot \vv$ does not add any new intersections, while
the intersection corresponding to $p$ is lost. This figure shows
$\ww=Babba$ (I) yielding
$\uu = aB$ and $\vv=bba$ (II).}\label{cut-the-curve}     
\end{figure}

Consider a minimal representative  of $w$ 
 drawn
in the fundamental domain. It must have self-intersections;
let $p$ be one of them. Let $\ww = x_1x_2\dots x_L$, (where $x_i \in \{a, A, b, B\}$), be a linear repsesentative for $w$, and 
suppose that $x_ix_{i+1}$ and $x_jx_{j+1}$, with $i<j$, 
are the two segments intersecting at $p$, (see \figc{cut-the-curve}, where
$x_ix_{i+1} = Ba$ and $x_jx_{j+1} = ba$). 
Set $\uu =x_{j+1}\dots x_L x_1x_2\dots x_i$ and 
$\vv = x_{i+1}\dots x_j$.  (In case $i+1=j$,
$\vv$ is a single-letter word). The cyclic words $u$ and $v$ 
together contain all the segments of $w$, except that $x_ix_{i+1}$ and $x_jx_{j+1}$ have
been replaced by $x_ix_{j+1}$ and $x_jx_{i+1}$.

Furthermore, there is a one-to-one correspondence between  the intersection points on $x_ix_{j+1}\cup x_jx_{i+1}$  and some subset of the intersection points on $x_ix_{i+1} \cup x_jx_{j+1}$. In fact, labeling  the
endpoints of the segment corresponding to  $x_ix_{i+1}$ (resp. $x_jx_{j+1}$) as $Q_i$ and $q_{i+1}$
(resp. $Q_j$ and $q_{j+1}$), as in \figc{cut-the-curve},
observe that the segment corresponding to  $x_ix_{j+1}$ and the broken arc $Q_i p q_{j+1}$ have the same endpoints, so any segment intersecting the first must intersect the second and therefore intersect part of $x_ix_{i+1} \cup x_jx_{j+1}$; similarly for $x_jx_{i+1}$ and
$Q_j p q_{i+1}$ (compare \figc{cut-the-curve}).
Therefore the change from $w$ to $u \cup v$ does not add any new intersections, while
the intersection corresponding to $p$ is lost. Hence
$\si(u)+\si(v)+1\leq \si(w).$ \end{proof}

The next lemma is needed in the proof of \propc{alpha}.

\start{lem}{six simple curves} In the doubly punctured plane $P$, 
if a reduced, non-empty word has a simple representative curve, then that 
curve is parallel to a boundary component. Thus with
the notation of \figc{pantsgrids} the only such
 words are ${a,b,ab,A,B \mbox{~and~} AB}$.
\end{lem}
\begin{proof} Let $\gamma$ be a simple, essential curve in $P$. 
Since $P$ is planar, $P \setminus \gamma$ has two connected components, 
$P_{1}$ and $P_{2}$. Since $\gamma$ is essential, neither $P_{1}$ nor $P_{2}$
is contractible, hence their Euler characteristics satisfy 
$\chi(P_{1}) \leq 0$ and $\chi(P_{2}) \leq 0$; since $\chi(P) = -1$
and $\chi(P)=\chi(P_{1})+\chi(P_{2})$ it follows that  either 
$\chi(P_{1})=0$ or $\chi(P_{2})=0$. Hence, one of the two connected 
components is an annulus, which implies that $\gamma$ is parallel to a boundary component, as desired.
\end{proof}

\start{prop}{alpha}  If $w$ is a positive 
cyclic
word representing a free homotopy class in the doubly punctured plane then $\si(w)  \geq \alpha(w)-1$ and  $\si(w)  \geq \beta(w)-1.$
\end{prop}
\begin{proof} 
By \lemc{six simple curves} the only words corresponding to simple curves 
are $a, b, ab$ and their inverses; for these, the statement
holds. In particular it holds for all words of length one. 
Suppose $w$ is any other positive word; it has length $L$ strictly
greater than 1. We may suppose by induction that the statement
holds for all words of length less than $L$.
By \propc{cut-the-word}, since the curve associated to $w$ is non-simple, 
the word $w$ has  
a linear representative $\ww$ which can be split as $\uu\cdot \vv$ 
so that the associated cyclic words satisfy $\si(w)\geq \si(u)+\si(v)+1$.
Note that $u$ and $v$ have length strictly less than $L$;
furthermore since $w$ is positive, so are $u$ and $v$.
Therefore by the induction hypothesis
$\si(u)+\si(v)+1\geq \alpha(u) -1 +\alpha(v) -1 +1$,  and
so $\si(w)  \geq \alpha(u) +\alpha(v) -1  = \alpha(w)-1.$
 The $\beta$ inequality is proved in the same way.
\end{proof}

\begin{prooftext}{Proof of \theoc{lower-bound}} By \propc{lower-case} 
there is a positive word $w'$ of length $L$ such that 
$\alpha(w') = \alpha(w), ~\beta(w') = \beta(w)$ and 
$\si(w) \ge \si(w')$. Then  
\propc{alpha} yields $\si(w') \geq \max \{\alpha(w),\beta(w)\}-1.$ 
Since $\alpha(w)+\beta(w)=L$ it follows that $\si(w) \geq L/2-1$ if $L$ is even and $\si(w) \geq (L+1)/2-1=(L-1)/2$ if $L$ is odd.
\end{prooftext}

\bibliographystyle{amsalpha}

{\sc

\noindent
Department of Mathematics,
Stony Brook University,
Stony Brook NY 11794
}

\noindent
\emph{E-mail}{:\;\;}\texttt{moira@math.sunysb.edu, tony@math.sunysb.edu}

\end{document}